\documentclass[12pt]{article}
\usepackage{amsmath, amssymb, bbm}

\newtheorem{theorem}{Theorem}[section]
\newtheorem{lemma}[theorem]{Lemma}
\newtheorem{corollary}[theorem]{Corollary}
\newtheorem{remark}[theorem]{Remark}
\newtheorem{proposition}[theorem]{Proposition}
\newtheorem{assumption}[theorem]{Assumption}

\newtheorem{definition}[theorem]{Definition}

\numberwithin{equation}{section}

\def\display{\displaystyle}
\def\C{{\mathbb C}}
\def\R{{\mathbb R}}

\def\E{{{\mathbb E}\,}}
\def\P{{\mathbb P}}

\def\eps{{\epsilon}}
\def\vareps{{\varepsilon}}
\def\phi{{\varphi}}
\def\lam{{\lambda}}
\def\th{{\theta}}
\def\al{{\alpha}}
\def\grad{{\nabla}}
\def\to{{\rightarrow}}

\def\proof{{\medskip\noindent {\bf Proof. }}}

\def\qed{{\hfill $\square$ \bigskip}}
\def\ms{\medskip}
\def\bs{\bigskip}
\def\q {\quad}
\def\qq {\qquad}
\def\nn{{\nonumber}}

  \def\sF {{\cal F}}
  
\def\sJ {{\cal J}}  \def\sL {{\cal L}}
\def\sM {{\cal M}}  
  \def\sR {{\cal R}}
\def\sS {{\cal S}}

  \def\bC {{\mathbb C}}
 \def\bE {{\mathbb E}}

  \def\bR {{\mathbb R}}

\begin{document}

\title{Uniqueness for the martingale problem associated with pure
jump processes of variable order}

\author{Huili Tang}

\maketitle

\begin{abstract}
\noindent Let $\sL$ be the operator defined on $C^2$ functions by
\[
\sL f(x)=\int[f(x+h)-f(x)-1_{(|h|\leq 1)}\grad f(x)\cdot
h]\frac{n(x,h)}{|h|^{d+\al(x)}}dh.
\]
This is an operator of variable order and the corresponding process
is of pure jump type.  We consider the martingale problem associated
with $\sL$. Sufficient conditions for existence and uniqueness are
given. Transition density estimates for $\al$-stable processes are
also obtained.

\end{abstract}

\vglue1.0truein {\small

\noindent{\it Subject Classifications:} Primary 60J75; Secondary
60G52

\noindent{\it Keywords:} martingale problem, symmetric stable
processes, jump process, transition density, homogeneous
distribution. }

\newpage

\section{Introduction}

We consider pure jump Markov processes corresponding to the
following infinitesimal generator:
\begin{equation}\label{defL}
\sL f(x)=\int[f(x+h)-f(x)-1_{(|h|\leq 1)}\grad f(x)\cdot
h]\frac{n(x,h)}{|h|^{d+\al(x)}}dh.
\end{equation}
The processes behave like a L\'{e}vy process at each point $x$, but
which process varies from point to point. Note that our operator
$\sL$ can be of variable order, i.e. $\alpha(x)$ is a function of
$x$.  In the case that either $n(x,h)$ or $\al(x)$ is a constant,
the corresponding process is called stable-like process. The
$1_{(|h|\leq 1)}\grad f(x)\cdot h$ term is omitted if $\al(x)<1$.
The question considered in this paper is the following. Is there a
process corresponding to the operator $\sL$, and if so, is there a
unique process\,?

In order to answer these questions, we consider the martingale
problem associated with $\sL$.
Let $\Omega=D[0,\infty)$ be the space of paths that are right
continuous with left limits, and let $\textit{X}_{t}:
\Omega\rightarrow \bR$ be defined by
$\textit{X}_{t}(\omega)=\omega(t)$. Let $\sF_{t}$ be the smallest
right continuous $\sigma$-field containing $\sigma(\textit{X}_{s},
s\leq t)$.
\ms  We say a probability measure $\P$ solves the martingale
problem for $\sL$ starting at $x_0$ if

\begin{description}
\item(i) $\P(X_0=x_0)=1,$ and
\item(ii) for every $f\in
C_{b}^{2}$, $f(X_t)-f(X_0)-\display\int_{0}^{t}\sL f(X_s)ds$ is a
$\P$-local martingale.
\end{description}

The purpose of this paper is to give sufficient conditions for the
existence and uniqueness of  the solution to the martingale problem
for pure jump processes of variable order.

There are only a very few papers \cite{Ba2}, \cite{Kolo},
\cite{Koma}, \cite{Nego}, \cite{Ts}, \cite{Ue} that handle variable
order terms without assuming a considerable amount of continuity in
the $x$ variable. Among the first is Bass \cite{Ba3}. The
infinitesimal generators of the processes considered there are given
by

\begin{align}\label{1-dim-gen}
A f(x)=\int[f(x+h)-f(x)-f^{'}(x)h1_{([-1,1])}(h)]v(x,dh).
\end{align}

As the result of a series of works \cite{Ba1}, \cite{Ba2},
\cite{Ba3} on stable-like processes in the $1980$'s, Bass \cite{Ba3}
handle the martingale problem for pure jump processes with variable
order using the Fourier transform.

Interested readers can find more details about jump processes in
\cite{Ba5}.

In Bass \cite{Ba3} a condition is given for uniqueness, and it is
stated there that there is no great difficulty in extending this to
higher dimensions. Unfortunately, the condition given is in terms of
the second derivative of the ratio of Fourier transforms, and can
really only be applied in the case where the jump kernel is of the
form $c|h|^{-1-\alpha(x)}\, dh$ for a suitable function $\alpha(x)$.
In recent years there has been considerable interest in operators
whose jump kernel is of the form $n(x,h)/|h|^{d+\alpha(x)}\, dh$,
where $n$ is a function that is bounded above and below; see
\cite{BBCK}, \cite{BL02}, \cite{BT08}, among others. For this reason
it is desirable to give a criterion for uniqueness directly in terms
of the functions $n$ and $\alpha(x)$, and that is the main purpose
of this paper.

The cases we consider in this paper are for multidimensional
processes and are much more general than \cite{Ba3}. We do a
perturbation of a multidimensional stable-like process. The
difficulty in this approach is threefold. The first difficulty is
that we have to establish new estimates on the transition densities
of multidimensional symmetric stable processes. Secondly, the
multidimensional case is much more singular than the one-dimensional
case. Lastly, the Fourier transform is hard to work with in our
case. Fortunately, we are mostly able to avoid the use of the
Fourier transform.

There are two perturbations involved in our proof. We first view
$\sL$ as a perturbation of stable-like processes, then we treat
stable-like processes as a perturbation of stable processes.

The rest of the paper is organized as follows. Section $2$ contains
notation, definitions, and statement of results. Section $3$
contains estimates on transition densities of $\al$-stable
processes. Some key estimates are obtained in Section $4$. Section
$5$ consists of the proof of uniqueness.

\section{Preliminaries}

We use the letter $c$ with subscripts to denote finite positive
constants whose exact values are unimportant and may change from
line to line. We use $C_{b}$ to denote the space of bounded
continuous functions on $\R^{d}$, and $C_{b}^{2}$ to denote the
space of bounded continuous functions on $\R^{d}$ that have bounded
continuous derivatives up to second order. The notation $C^{\infty}$
denote the collection of functions that have continuous derivatives
of any order. Let $C_{c}$ denote the space of functions in $C_{b}$
with compact support and similarly let $C_{c}^{\infty}$ and
$C_{c}^{2}$ denote the collections of functions in $C^{\infty}$ and
$C^{2}$, respectively, with compact support. We denote
$\displaystyle\sup_{x} |f(x)|$ by $\|f\|$. We use $|x|$ to denote
the Euclidean norm for $x \in \R^d$. The notation $:=$ is to be read
as "is defined to be." For two real numbers $a$ and $b$, $a\land b
:=\min\{a, b\}$. For a function $f$ on $\R^{d}$, its Fourier
transform $\hat{f}$ is defined by
\[
\hat{f}(u) := \int_{\R^{d}}e^{iu\cdot x}f(x)dx, \qq u\in \R^{d}.
\]

A multidimensional symmetric stable process of index $\al$ is a
L\'{e}vy process $X_t$ such that
\[
\E e^{iu\cdot Z_{t}}=e^{-t|u|^{\al}}.
\]
The L\'{e}vy measure for such a process is given by
$\frac{c_{\al}}{|h|^{d+\al}}dh$, where $c_{\al}$ is a constant
depending only on $\al$. This follows because the
L\'{e}vy-Khintchine formula says that
\[
\E e^{iu\cdot Z_{t}}=e^{-t\Phi(u)},
\]
where
\[
\Phi(u)=\int_{|h|\neq 0}(e^{iu\cdot h}-1-iu\cdot
h1_{(|h|\leq1)})\frac{c_{\al}}{|h|^{d+\al}}dh.
\]
With a change of variables $h=\frac{v}{|u|}$ and the fact
$h/|h|^{d+\al}$ is odd, we have
\[
\Phi(u)=|u|^{\al}\int_{|v|\neq 0}(e^{i\frac{u}{|u|}\cdot
v}-1-i\frac{u}{|u|}\cdot
v1_{(|v|\leq|u|)})\frac{c_{\al}}{|v|^{d+\al}}dv=c|u|^\al.
\]

For the existence of a solution to the martingale problem, we need
the following assumptions.
\begin{assumption}\label{existassum}
Suppose
\begin{description}
\item{(a)} for all $x$, there exist positive constants $c_1$, $c_2$
such that $c_1\leq n(x,h)\leq c_2$.
\item{(b)} $\sL f$ is continuous whenever $f\in C^2_b$.
\end{description}
\end{assumption}

For the uniqueness of the solution to the martingale problem for
$\sL$ as defined in (\ref{defL}), we need the following assumption.

\begin{assumption} \label{uniqassum}
Suppose
\begin{description}
\item{(a)}  there exist positive constants $c_1$, $\gamma$ and $\eps$ and
a Dini continuous function $\xi:\R^d\to (0,\infty)$ such that for
all $x$, $|n(x,h)-\xi(x)| \leq c_1(1\land |h|^{\eps}) $;
\item{(b)} $0<\underline{\al}=\display\inf_{x}\al(x)\leq\display\sup_{x}\al(x)=\overline{\al}<2$;
\item{(c)} $\beta(z)=o(1/|\ln z|)$ as $z \to 0,$ where $\beta(z)=\display\sup_{|x-y|\leq
z}|\alpha(x)-\alpha(y)|$;
\item{(d)} $\int^{1}_{0} \frac{\beta (z)}{z^{1+\gamma}}< \infty$.
\end{description}
\end{assumption}

We say a function $\xi(x)$ is Dini continuous if
\[
\int_{0}^{1}\frac{\psi(z)}{z}dz<\infty,\q \mathrm{where} \,\,
\psi(z)=\display\sup_{|x-y|\leq z}|\xi(x)-\xi(y)|.
\]

We also temporarily assume the following on $\xi(x)$.

\begin{assumption}\label{tempassum}
There exists a positive constant $\zeta$ such that
\[
|\xi(x)|\leq \zeta, \qq x \in \R^d.
\]
\end{assumption}

Our existence theorem is the following.

\begin{theorem}\label{existence}
Suppose that Assumption \ref{existassum} holds. Then for every $x_0
\in \bR$ there exists a solution to the martingale problem for $\sL$
starting from $x_0$.
\end{theorem}

\proof Bass \cite{Ba3} gives a complete proof for the existence for
one-dimensional case and it has no difficulty to extend the same
proof to higher dimensions. The idea in the proof is to construct a
sequence of tight probability measures $\P_{n}$ and show there is a
subsequence of $\P_{n}$ which converges to $\P$, a solution to the
martingale problem.

\qed

Our main result for uniqueness is the following.

\begin{theorem}\label{main thm}
Suppose that Assumption \ref{uniqassum} holds. Then for each $x_0$
the martingale problem associated with the operator $\sL$ starting
at $x_0$ has a unique solution.
\end{theorem}

The conditions in our uniqueness theorem are quite mild. A recent
paper by Barlow et al \cite{BBCK} indicates that
uniqueness can fail if one only requires that  $n(x,h)$ be bounded.

\section{Transition densities of $\al$-stable processes}

In this section, we will obtain a power series expansion for the
transition density of a symmetric stable process in $d$ dimensions.

The estimate (\ref{estp1}) on the transition density of a symmetric
stable process is known; see Kolokoltsov \cite[Proposition
3.1]{Kolo}. But we prove it using a different approach. Our approach
allows us to obtain an estimate on the second derivative of
transition density $p_{t}(x,y)$ by differentiating the power series.

Let $0<\al<2$ be fixed, let $X_t$ be a multidimensional symmetric
$\al$-stable process, and let $p_{t}(x,y)$ be the transition density
of $X_t$. The characteristic function of $X_1$ is $\E
\exp(iuX_1)=\exp(-|u|^{\al}).$ Let $u=(u_1,u_2,...,u_d)$ be a vector
and $\beta=(\beta_1,\beta_2,...,\beta_d)$ be a multi-index with
nonnegative integers entries; define the size of a multi-index
$\beta$ by $|\beta|=|\beta_1|+...+|\beta_d|$ and define
$u^{\beta}=\prod_{j=1}^{d}u_{j}^{\beta_j}$ and
$\partial^{\beta}f=\partial_{1}^{\beta_{1}}...\partial_{d}^{\beta_{d}}f$.
Since $u^{\beta}\exp(-|u|^{\al})$ is integrable for all
multi-indices $\beta$, $p_{1}(0,x)$ has bounded partial derivatives
of all orders. We have the following estimates. The proof partially
follows \cite{Grafakos}.

\begin{proposition}\label{tranden}
There exist positive constants $c_1$ and $M_1$ such that if $|x|\ge
M_1$,
\begin{align}
p_{1}(0,x)=(2
\pi)^{-d}\sum_{k=1}^{\infty}\frac{c_{d,k\al}}{k!}|x|^{-(d+k\al)},\label{dpe}
\end{align}
where $c_{d,k\al}=2^{k\al}\pi^{-d/2}
\frac{\Gamma((d+k\al)/2)}{\Gamma(-k\al/2)}$. Furthermore
\begin{align}
|p_{1}(0,x)|=c_1 |x|^{-(d+\al)}(1+o(1)).\label{estp1}
\end{align}
\end{proposition}

Before we proceed to the the proof, we give a definition of a
homogeneous distribution which is needed in our proof.
\begin{definition}
Suppose $f$ is in the Schwartz class. For $z \in \bC$, a
homogeneous distribution $u_{z}$ is defined as follows:
\begin{align}
u_{z}(f)=\int_{\bR^{d}}\frac{\pi^{\frac{z+d}{2}}}{\Gamma(\frac{z+d}{2})}|x|^{z}f(x)dx.\label{homodistr}
\end{align}
\end{definition}

It is clear that the integral converges for $Re~z>-d$. We would like
to extend the definition of $u_{z}(f)$ to all $z \in \bC$. Let
$Re~z>-d$ and $N$ be a fixed positive integer. For $f$ in the
Schwartz class, rewrite the integral in (\ref{homodistr}) as
follows.

\begin{align}
&~~\int_{|x|<1}\frac{\pi^{\frac{z+d}{2}}}{\Gamma(\frac{z+d}{2})}|x|^{z}
\Big\{f(x)-\sum_{|\beta| \leq N}\frac{(\partial^{\beta} f)(0)}{\beta ! }x^{\beta} \Big\}dx \label{homodistr-1}\\
&\qq+\int_{|x|<
1}\frac{\pi^{\frac{z+d}{2}}}{\Gamma(\frac{z+d}{2})}|x|^{z}
\sum_{|\beta|\leq N}\frac{(\partial^{\beta} f)(0)}{\beta !} x^{\beta} dx \label{homodistr-2}\\
&\qq+\int_{|x|\ge
1}\frac{\pi^{\frac{z+d}{2}}}{\Gamma(\frac{z+d}{2})}|x|^{z}f(x)dx.
\label{homodistr-3}
\end{align}

Suppose $Re~z>-N-d-1$. Since the difference inside the brackets of
(\ref{homodistr-1}) is bounded by a constant multiple of
$|x|^{N+1}$, the integral in (\ref{homodistr-1}) is a well defined
analytic function. It is obvious that the integral in
(\ref{homodistr-3}) is also well defined since $f$ is in the
Schwartz class. For the integral in (\ref{homodistr-2}), we use
polar coordinates to get
\begin{align}
&\frac{\pi^{\frac{z+d}{2}}}{\Gamma(\frac{z+d}{2})}
\sum_{|\beta|\leq N}\frac{(\partial^{\beta} f)(0)}{\beta !}
\int_{0}^{1}\int_{\sS}(r\th)^{\beta}r^{z+d-1}dr\,d\th, \nn \\
&=\frac{\pi^{\frac{z+d}{2}}(\partial^{\beta}
f)(0)}{\Gamma(\frac{z+d}{2})}\frac{1}{\beta !}(\int_{\sS}
\th^{\beta}d\th)\int_{0}^{1}r^{|\beta|+z+d-1}dr, \nn \\
&=\frac{\pi^{\frac{z+d}{2}}(\partial^{\beta} f)(0)}{\beta !
\Gamma(\frac{z+d}{2})}
\frac{\int_{\sS}\th^{\beta}d\th}{|\beta|+z+d}.\label{homodistr-2-polar}
\end{align}
The integral in (\ref{homodistr-2-polar}) is zero when $|\beta|$ is
odd. If $|\beta|$ is even, $(~|\beta|+z+d)^{-1}$ has a simple pole
at $z=-d-|\beta|$ for $|\beta|\leq N $ and even. We know that
$\Gamma(\frac{z+d}{2})$ also has a simple pole at
$z=-d-2j,~j=1,2,...[\frac{N}{2}]$. These poles exactly cancel with each
other. We therefore see that the integral in
(\ref{homodistr}) is a well defined analytic function when
$Re~z>-N-d-1$. Since N was arbitrary, (\ref{homodistr}) is well
defined for all $z \in \bC$.

We have the following lemma.
\begin{lemma}\label{iftlemma}
For all positive integers $k$
\[
\lim_{\eps \to 0}\int_{\R^d}e^{-iux}|u|^{k\al}e^{-\eps
|u|^{2}/2}du=(2\pi)^{-d}\frac{c_{d,k\al}}{|x|^{(d+k\al)}}.
\]
\end{lemma}
\proof First we look at a more general case of the above. For
all $z\in \C$ and $f$ in the Schwartz class, we will use polar
coordinates $x=r\th$ and $u=t\phi$. The following is justified by
Fubini and rotational invariance.
\begin{align*}
&~\int_{\R^d}|u|^{z}\hat{f}(u)du=\int_{\R^d}\int_{\R^d}|u|^{z}e^{ix\cdot u}f(x)\,dx\,du\\
&=\int_{0}^{\infty}\int_{\sS}\int_{0}^{\infty}\int_{\sS}e^{irt\th\cdot\phi}\,d\phi \,t^{d+z-1}dt\,f(r\th)d\th\, r^{d-1}\,dr\\
&=\int_{0}^{\infty}r^{-(d+z)}\int_{\sS}\int_{0}^{\infty}\int_{\sS}e^{it\th\cdot\phi}\,d\phi \,t^{d+z-1}\,dtf(r\th)\,d\th \,r^{d-1}\,dr\\
&=\int_{0}^{\infty}r^{-(d+z)}\int_{\sS}\int_{0}^{\infty}\int_{\sS}e^{it\phi_{1}}\,d\phi \,t^{d+z-1}dt\,f(r\th)d\th \,r^{d-1}dr\\
&=c_{d,z}\int_{0}^{\infty}r^{-(d+z)}\int_{\sS}f(r\th)d\th \,r^{d-1}dr\\
&=c_{d,z}\int_{\R^d}|x|^{-(d+z)}f(x)dx,
\end{align*}
where $\phi_1$ is the first coordinate of $\phi$,
\begin{align}
\sigma(t)=\int_{\sS}e^{it\phi_{1}}d\phi\label{sigmat}.
\end{align}
and
\begin{align}
c_{d,z}=\int_{0}^{\infty}\int_{\sS}e^{it\phi_{1}}d\phi\,
t^{d+z-1}dt=\int_{0}^{\infty}\sigma(t)t^{d+z-1}dt.\label{c(d,z)}
\end{align}

Next we need to  show that $c_{d,z}$ is bounded for some range of
$z$'s. After doing a change of variable, we get
\[
\sigma(t)=\int_{-1}^{1}e^{its}\omega_{d-2}(\sqrt{1-s^2})^{d-2}\frac{ds}{\sqrt{1-s^2}}=c_{d}J_{d-2}(t).
\]
Using the asymptotics for Bessel functions, we get that
$|\sigma(t)|\leq c t^{-1/2}$ when $d-2>-1/2$. If $-d<Re~z<-d+1/2$,
then
\begin{align*}|
c_{d,z}|&\leq\int_{0}^{\infty}|\sigma(t)|t^{Re~z+d-1}dt\\
&\leq\int_{0}^{1}\omega_{d-1}t^{Re~z+d-1}dt+c_{d}\int_{1}^{\infty}t^{Re~z+d-3/2}dt\\
&<\infty.
\end{align*}
Since the function
$z\to\int_{\R^d}|u|^{z}\hat{f}(u)du-c_{d,z}\int_{\R^d}|x|^{-(d+\al)}f(x)dx$
is entire and vanishes for $-d<Re~z<-d+1/2$ and every $f$ in the Schwartz class, it must vanish everywhere.\\

Now letting $z=k\al$ for $k=1,2,...$ and
$f(y)=f(y-x)=e^{-|y-x|^2/2\eps}$, we obtain
\begin{align}
\int_{\R^{d}}e^{-iu\cdot
x}|u|^{k\al}e^{-\eps|u|^2/2}du=c_{d,k\al}(2\pi \eps)^{-d/2}
\int_{\R^d}|y|^{-(d+k\al)}e^{-|x-y|^2/2\eps}dy.\label{ift}
\end{align}
Letting $\eps\to 0$ in (\ref{ift}), we have
\begin{align}
\lim_{\eps \to 0}(2\pi)^{-d}\int_{\R^{d}}e^{-iu\cdot
x}|u|^{k\al}e^{-\eps|u|^2/2}du=(2\pi)^{-d}c_{d,k\al}
|x|^{-(d+\al)}.\label{1irt1}
\end{align}
\qed

\bs We now prove Proposition \ref{tranden}.

\proof By the Fourier inversion theorem, we know that\\
\begin{align*}
p_{1}(0,x)&=\Big(\frac{1}{2\pi}\Big)^{d}\int_{\R^d} e^{-iu\cdot x}e^{-|u|^{\al}}du\\
&=\lim_{\eps\to 0}(2\pi)^{-d}\int_{\R^d} e^{-iu\cdot
x}(e^{-|u|^{\al}}-1)e^{-\eps|u|^{2}/2}du\\
&~~~+\lim_{\eps\to 0}(2\pi)^{-d}\int_{\R^d} e^{-iu\cdot x}e^{-\eps|u|^{2}/2}du\\
&=I_1+I_2.
\end{align*}

Looking at $I_2$,
\begin{align}
I_2=\lim_{\eps\to 0}(2\pi\eps)^{-d/2}e^{-|x|^2/2\eps}=0 ~~~if~ x
\neq 0.\label{ift2}
\end{align}

Next, looking at $I_1$ and using the Taylor expansion of $e^{x}$, we
get
\begin{align}
I_1=\lim_{\eps\to
0}(2\pi)^{-d}\displaystyle\sum_{k=1}^{\infty}\int_{\R^d} e^{-iu\cdot
x}(-1)^{k} \frac{|u|^{k\al}}{k!}e^{-\eps|u|^{2}/2}du.\label{ift1}
\end{align}

Then applying Lemma \ref{iftlemma} in (\ref{ift1}), we have
\begin{align}\label{ift1-expansion}
I_1=(2\pi)^{-d}\sum_{k=1}^{\infty}(-1)^{k}\frac{c_{d,k\al}}{k!}|x|^{-(d+k\al)}|x|^{-(d+\al)}.
\end{align}

The first part of Proposition \ref{tranden} is now proved by
combining \eqref{ift2} and \eqref{ift1-expansion}.

It remains to show
\begin{align*}\label{series-sum}
\sum_{k=2}^{\infty}(-1)^{k}\frac{c_{d,k\al}}{k!}|x|^{-(d+k\al)}=O(|x|^{-(d+2\al)}).
\end{align*}

As we can see from Remark \ref{remark:1},
\[
c_{d,k\al}=2^{k\al}\pi^{-d/2}
\frac{\Gamma((d+k\al)/2)}{\Gamma(-k\al/2)}.
\]

For convenience, we set the series coefficients
$(-1)^{k}\frac{c_{d,k\al}}{k!}=a_{k}$ for any $k$.

By Stirling's formula, we have
\[
\lim_{x \to \infty} \frac{\Gamma(x+1)}{({x \over e})^x \sqrt{2\pi
x}}=1.
\]

Applying Stirling's formula and the fact that ${\al\over 2}<1 $, we
get

\begin{align*}
&~~~\lim_{k\to \infty}{a_{k+1} \over a_{k}}\\
&=\lim_{k \to \infty}
\frac{\Gamma({d+(k+1)\al \over 2})}{(k+1)\Gamma({d+k\al \over 2})}\\
&=\lim_{k \to \infty}\frac{1}{k+1}\Big(\frac{(d+k\al)/
2+\al/2-1}{(d+k\al)/ 2-1}\Big)^{\big(\tfrac{d+k\al}{
2}+\tfrac{\al}{2}-1\big)}\Big(\frac{(d+k\al)/
2+\al/2-1}{e}\Big)^{\al/2}\\
&=e^{\al/2}\lim_{k \to \infty}k^{\al/2-1}\\
&=0.
\end{align*}

This completes the proof of Proposition \ref{tranden}.\qed

An alternative approach is to view the symmetric stable process as
Brownian motion subordinated by a one-sided one-dimensional stable
process of index $\alpha/2$, and to use the known density for these
one-sides processes. Although there is an explicit expression
available for the latter, it is given as an infinite series, so this
method does not seem any shorter or simpler than ours.

\section{Estimates}

In this section, we
will obtain some key estimates, which will be
used in our proof of uniqueness.

\begin{proposition}\label{estd}
There exists a positive constant $c_{1}>0$ which depends on $d$ and
$\al$ such that
\begin{description}
\item{(a)} $|p_{1}(0,x)|\leq c_{1} (1\wedge |x|)^{-(d+\al)}$ and for $k=1,2,$
\[
| \partial^{k} p_{1}(0,x)| \leq c_{1} (1\wedge |x|^{-(d+\al+k)}).
\]
\item{(b)} $|p_{t}(0,x)|\leq c_{1} t^{-d/\al}\big(1\wedge (t^{1/\al}|x|^{-1})\big)^{(d+\al)}$ and for $k=1,2,$
\[
| \partial^{k} p_{t}(0,x)| \leq c_{1} t^{-(d+k)/\al}\big(1\wedge
(t^{1/\al}|x|^{-1})\big)^{(d+\al)+k}.
\]
\end{description}
\end{proposition}
\proof (b) follows from (a) by scaling. For (a), the first estimate is just a restatement of Proposition
\ref{tranden}. We have the full expansion in Proposition \ref{tranden}. Differentiating it with respect to $x$
and following a similar argument to proving (\ref{ift1-expansion}) gives the case $k=1,2$.
\qed

\ms

Fix $\lam>0$, and for bounded $f$ let
\begin{align}
\sR_{\lam} f(x)=\bE\int_{0}^{\infty}e^{{-\lambda}t}f(X_{t})dt=\int
r^{\lambda}(x-y)f(y)dy.\label{resolvent}
\end{align}
where $X_t$ is a symmetric stable process, $p_t(0,x)$ is its
transition probability, and $r^{\lam}(x)=\int_{0}^{\infty} e^{-\lam
t} p_t(0,x)dt$. We also let
\begin{align}
\sM_{z} f(x) = \int[f(x+h)-f(x)-1_{(|h|\leq 1)}\nabla f(x)\cdot
h]\frac{\xi (z)}{|h|^{d+\alpha(z)}}dh. \label{stablegenerator}
\end{align}
for $f\in C_{b}^{2}$. Observe that in (\ref{stablegenerator}), $z$
is a parameter and $\sM_{z}f$ is a function only of $x\in \bR^{d}$.

We will investigate the operator $\sM_{y}$, which is the generator
of a symmetric stable process of fixed order $\al(y)$. In this case,
the operator $\sM_{y}$ has L\'{e}vy measure $\frac{\xi
(y)}{|h|^{d+\al(y)}}dh$. We define $\sR_{\lam}^{c}$ and
$r_y^{\lam,c}$ by (\ref{resolvent}) when the process $X_t$ of
(\ref{resolvent}) is generated by $\sM_{y}$.

We have the following estimates regarding the resolvent density of
$\sR_{\lam}^{c}$, where for convenience we write $\al$ in place of
$\al(y)$.

\begin{proposition}\label{densityofRc}
There exists a positive constant $c_{0} \in (0,\infty)$
such that the following hold:
\begin{description}
\item{(a)} $r_y^{\lam,c}(x)\leq c_{0}(\tfrac{1}{\lam}|x|^{-2\al}\land 1)|x|^{-d+\al}$;
\item{(b)} $\sum_i|\partial r_y^{\lam,c}(x)/\partial x_i|\leq
c_{0}(\tfrac{1}{\lam}|x|^{-2\al}\land 1)|x|^{-d+\al-1}$;
\item{(c)} $\sum_{i,j}|\partial^{2}r_y^{\lam,c}(x)/\partial x_i \partial x_j|\leq
c_{0}(\tfrac{1}{\lam}|x|^{-2\al}\land 1)|x|^{-d+\al-2}$.
\end{description}
\end{proposition}

\proof We will only prove part (a), the others being similar. We
know by Proposition \ref{estd} that there exists a positive constant
$c_{1}$ such that $p_{t}(0,x)\leq c_{1} t^{-d/\al}(1 \wedge
(t^{1/\al})|x|^{-1})^{(d+\al)}$. Then
\begin{align*}
r_y^{\lam,c}(x)
&=\int_{0}^{\infty} e^{-\lam t}p_{t}(0,x)dt\\
&\leq c_{1} |x|^{-(d+\al)} \int_{0}^{|x|^{\al}}te^{-\lam t}\,dt+ c_{1}
\int_{|x|^{\al}}^{\infty}e^{-\lam
t}t^{-d/\al}dt\\
&\leq I_1 + I_2.
\end{align*}

First, we consider $|x|\geq 1$. For $I_1$,
\[
I_1 \leq c_{2} |x|^{-(d+\al)}\int_{0}^{|x|^{\al}}e^{-\lam t/2}dt
\leq c_{3}\lam^{-1} |x|^{-(d+\al)}.
\]

Next,  since $e^{-\lam |x|^{\al}}\leq c_4|x|^{-\al}$ when $|x|\geq 1$,
\[
I_2 \leq c_{5}|x|^{-d}\int_{|x|^{\al}}^{\infty}e^{-\lam t}dt \leq
c_{6}\lam^{-1} |x|^{-d}e^{-\lam |x|^{\al}}\leq
c_{7}\lam^{-1}|x|^{-(d+\al)}.
\]
Summing $I_1$ and $I_2$, we get for $|x|\geq 1$,
\[
r_y^{\lam,c}(x)\leq (c_3+c_7)\lam^{-1} |x|^{-(d+\al)}.
\]

A similar proof also works for $|x|\leq 1$.
Again look at $I_1$ and $I_2$.
\[
I_1 \leq c_{1} |x|^{-(d+\al)}\int_{0}^{|x|^{\al}}te^{-\lam t}\,dt
\leq c_{8}|x|^{-d+\al},
\]
and
\[
I_2\leq c_{1} \int_{|x|^{\al}}^{\infty}t^{-d/\al}e^{-\lam t}dt\leq
 c_9 e^{-\lam |x|^\al}\int_{|x|^\al}^\infty t^{-d/\al}\, dt
=c_{10} |x|^{-d+\al}.
\]

Summing $I_1$ and $I_2$, we get for $|x|<1$,
\[
r_y^{\lam,c}(x) \leq (c_8+c_{10}) |x|^{-d+\al}.
\]
The two cases above prove the estimates. \qed

Let $\varphi$ be an even radial nonnegative $C^\infty$ function with
support in $B(0,1/2)$ and $\int\varphi(x)dx=1$. Define
$\varphi_{\vareps}=\vareps^{-d}\varphi(x/\vareps)$. Let
$\lam\in[1,\infty)$ be fixed. Define:

\begin{align}
r^{\lam,\varepsilon}_{y}=\E^{x}\int^{\infty}_{0}e^{-{\lam}
t}\varphi_{\varepsilon}(\textit{X}^{y}_{t})dt,
\end{align}
where $\textit{X}^{y}_{t}$ is a stable process generated by
(\ref{stablegenerator}) with L\'evy measure
$\frac{c}{|h|^{d+\al(y)}}dh$
Then we have the following estimates on $r^{\lam,\varepsilon}_{y}$.

\begin{proposition}\label{est-on-r^eps}
There exists a positive constant $c_{0} \in (0,\infty)$ such that
the following hold:
\begin{description}
\item{(a)} $r^{\lam,\varepsilon}_{y}(x)\leq c_{0}(\tfrac{1}{\lam}|x|^{-2\al}\land 1)|x|^{-d+\al}$;
\item{(b)} $\sum_i|\partial r^{\lam,\varepsilon}_{y}(x)/\partial x_i|\leq
c_{0}(\tfrac{1}{\lam}|x|^{-2\al}\land 1)|x|^{-d+\al-1}$;
\item{(c)} $\sum_{i,j}|\partial^{2}r^{\lam,\varepsilon}_{y}(x)/\partial x_i \partial x_j|\leq
c_{0}(\tfrac{1}{\lam}|x|^{-2\al}\land 1)|x|^{-d+\al-2}$.
\end{description}
\end{proposition}

\proof Again, we will only prove part (b) as the others are similar.
To get estimates on $r^{\lam,\varepsilon}_{y}(x)$, we write
\begin{align*}
|r^{\lam,\varepsilon}_{y}(x)|&=\Big | \int
r^{\lam,c}_{y}(x-u)\varphi_{\vareps}(u)du \Big|\\
&\leq \Big| \int
\big(r^{\lam,c}_{y}(x-u)-r^{\lam,c}_{y}(x)\big)\varphi_{\vareps}(u)du
\Big| +\Big|\int r^{\lam,c}_{y}(x)\varphi_{\vareps}(u)du \Big|\\
&=I_1+I_2.
\end{align*}

We estimate $I_1$ first.
\begin{align*}
I_1& \leq \int_{|u|\leq \frac{|x|}{2}}
\big|\big(r^{\lam,c}_{y}(x-u)-r^{\lam,c}_{y}(x)\big)\varphi_{\vareps}(u)\big|du\\
&~~~+\int_{\frac{|u|}{2}<|u|\leq \frac{3|x|}{2}}
\big|\big(r^{\lam,c}_{y}(x-u)-r^{\lam,c}_{y}(x)\big)\varphi_{\vareps}(u)\big|du\\
&~~~+\int_{|u|>\frac{3|x|}{2}}
\big|\big(r^{\lam,c}_{y}(x-u)-r^{\lam,c}_{y}(x)\big)\varphi_{\vareps}(u)\big|du\\
&=I_{11}+I_{12}+I_{13}.
\end{align*}

Since $\displaystyle\sup_{|u| \leq
\frac{|x|}{2}}|r^{\lam,c}_{y}(x-u)-r^{\lam,c}_{y}(x)|\leq
c_{1}|u|\sum_i|\partial r^{\lam,\varepsilon}_{y}(x/2)/\partial
x_i|$, we have
\[
I_{11} \leq c_{1}\sum_i|\partial
r^{\lam,\varepsilon}_{y}(x/2)/\partial x_i|\int_{|u|\leq
\frac{|x|}{2}}|u|\varphi_{\vareps}(u)du \leq
c_{2}|r^{\lam,\varepsilon}_{y}(x/2)|.
\]

As for $I_{12}$, since $\varphi_{\vareps}(x)$ has support
$B(0,1/2\vareps)$, we have

\begin{align*}
I_{12}&\leq \int_{\frac{|x|}{2}<|u|\leq \frac{3|x|}{2}}\big\{
|r^{\lam,c}_{y}(x-u)|+|r^{\lam,c}_{y}(x)|\big\}\varphi_{\vareps}(u)du\\
&\leq c_{3}|r^{\lam,c}_{y}(x)|.
\end{align*}

Looking at $I_{13}$, since $|x-u|>|x|/2$ when $|u|>3|x|/2$, we have
\[
I_{13}\leq c_{4}
|r^{\lam,c}_{y}(x)|\int_{|u|>\frac{3|x|}{2}}\varphi_{\vareps}(u)du\leq
c_{5}|r^{\lam,c}_{y}(x)|.
\]

It is easy to see that $I_2=r^{\lam,c}_{y}$ since
$\int\varphi(x)dx=1$.
Using Proposition \ref{densityofRc} and combining with the estimates for
$I_1$ and $I_2$
finishes the proof.

\qed

\begin{corollary}\label{est-on-r^eps-coro}
There exists a positive constant $c_{0} \in (0,\infty)$ such that
the following hold:
\begin{description}
\item{(a)} $r^{\lam,\varepsilon}_{y}(x)\leq c_0 \tfrac{1}{\lam}|x|^{-d-\al}$;
\item{(b)} $\sum_i|\partial r^{\lam,\varepsilon}_{y}(x)/\partial x_i|\leq c_0 \tfrac{1}{\lam}|x|^{-d-\al-1}$;
\item{(c)} $\sum_{i,j}|\partial^{2}r^{\lam,\varepsilon}_{y}(x)/\partial x_i \partial x_j|\leq
c_0\tfrac{1}{\lam}|x|^{-d-\al-2}$.
\end{description}
\end{corollary}

\proof The Corollary follows easily by looking at small $|x|$ in
Proposition \ref{est-on-r^eps}. \qed

We now have the following proposition.
\begin{proposition}\label{pertur-on-lc}
If Assumption \ref{uniqassum} (a) holds, there exist positive
constant $\eta$ and $k_{1} \in (0.\infty)$ such that
\begin{eqnarray*}
  |(\sL-\sM_{x})r^{\lam,\varepsilon}_{y}(u)| &\leq & k_1 \frac{1}{|u|^{d+\al(x)-\al(y)-\eta}},\qq |u|\leq 1,\\
  |(\sL-\sM_{x})r^{\lam,\varepsilon}_{y}(u)| &\leq & \frac{k_1}{\lam}\frac{1}{|u|^{d+\al(x)+\al(y)}},\qq |u|>
  1.
\end{eqnarray*}

In particular, for all $u$
\[
|(\sL-\sM_{x})r^{\lam,\varepsilon}_{y}(u)|\leq
\frac{k_1}{\lam}\frac{1}{|u|^{d+\al(x)+\al(y)}}.
\]
\end{proposition}
\proof For convenience, we set
\begin{align}\label{J(u,h)}
J(u,h)=r_y^{\lam,\varepsilon}(u+h)-r_y^{\lam,\varepsilon}(u)-\grad
r_y^{\lam,\varepsilon}(u)
  \cdot h 1_{(|h|\leq 1)}
\end{align}
When $|u|\leq 1$, we have
\begin{align*}
  &~~~~|(\sL-\sM_x)r^{\lambda,\varepsilon}_{y}(u)|\\
  &=\Big|\int J(u,h) \frac{n(x,h)-\xi(x)}{|h|^{d+\al(x)}}dh\Big|\\
  & \leq \int_{|h|\leq \frac{|u|}{2}}\Big|J(u,h)\frac{n(x,h)-\xi(x)}{|h|^{d+\al(x)}}\Big|dh
  +\int_{\frac{|u|}{2}< |h| \leq \frac{3|u|}{2}}\Big|J(u,h)\frac{n(x,h)-\xi(x)}{|h|^{d+\al(x)}}\Big|dh\\
  &~~~+\int_{|h|>\frac{3|u|}{2}}\Big|J(u,h)\frac{n(x,h)-\xi(x)}{|h|^{d+\al(x)}}\Big|dh\\
  &=I_1+I_2+I_3.
\end{align*}

First of all, looking at $I_1$, by Assumption \ref{uniqassum} (a)
and Proposition \ref{est-on-r^eps} we have
\[
I_1 \leq c_{2}\sup_{
B(u,|u|/2)}|\partial^{2}r^{\lam,\varepsilon}_{y}(z)|\int_{|h|\leq
\frac{|u|}{2}} \frac{1}{|h|^{d+\al(x)-2-\eps}}dh \leq
\frac{c_{3}}{|u|^{d+\al(x)-\al(y)-\eps}}.
\]
Next for $I_2$, by Proposition \ref{est-on-r^eps} we have
\[
|J(u,h)|\leq
\frac{c_{4}}{|u|^{d-\al(y)}}+\frac{c_{5}}{|u+h|^{d-\al(y)}}.
\]

Thus we have
\begin{align*}
I_2 & \leq \int_{\frac{|u|}{2}< |h| \leq
\frac{3|u|}{2}}\Big\{\frac{c_{4}}{|u|^{d-\al(y)}}+\frac{c_{5}}{|u+h|^{d-\al(y)}}\Big\}
\frac{|n(x,h)-\xi(x)|}{|h|^{d+\al(x)}}dh\\
&\leq\frac{c_6}{|u|^{d+\al(x)-\al(y)-\eps}}+\frac{c_{7}}{|u|^{d+\al(x)-\eps}}
\int_{\frac{|u|}{2}< |h| \leq
\frac{3|u|}{2}}\frac{1}{|u+h|^{d-\al(y)}}dh\\
&\leq \frac{c_{8}}{|u|^{d+\al(x)-\al(y)-\eps}}.
\end{align*}

For $I_3$, there are two cases.\\

Case 1: If $\frac{3|u|}{2}< 1$, we break $I_3$ into two pieces as
follows:
\begin{align*}
I_3&=\int_{\frac{3|u|}{2}<|h|\leq
1}|J(u,h)|\frac{|n(x,h)-\xi(x)|}{|h|^{d+\al(x)}}dh+\int_{|h|>
1}|J(u,h)|\frac{|n(x,h)-\xi(x)|}{|h|^{d+\al(x)}}dh\\
&=I_{31}+I_{32}.
\end{align*}
We assume $\al(x)\geq 1$, then
\begin{align*}
I_{31} &\leq \int_{\frac{3|u|}{2}<|h|\leq
1}\Big\{\frac{c_{9}}{|u+h|^{d-\al(y)}}+\frac{c_{10}}{|u|^{d-\al(y)}}+\frac{c_{11}|h|}{|u|^{d-\al(y)+1}}\Big\}
|h|^{-d-\al(x)+\eps}dh\\
& \leq
\frac{c_{12}}{|u|^{d+\al(x)-\al(y)-\eps}}+\frac{c_{13}}{|u|^{d-\al(y)+1}}
\int_{\frac{3|u|}{2}<|h|\leq 1}|h|^{-d-\al(x)+1+\eps}dh\\
&\leq
\frac{c_{12}}{|u|^{d+\al(x)-\al(y)-\eps/2}}+\frac{c_{13}}{|u|^{d+\al(x)-\al(y)-\eps/2}}
\int_{\frac{3|u|}{2}<|h|\leq 1} |h|^{-d+\eps/2}dh\\
&\leq \frac{c_{14}}{|u|^{d+\al(x)-\al(y)-\eps/2}},
\end{align*}
since $|h|^{-\al(x)+1+\eps/2}\leq |u|^{-\al(x)+1+\eps/2} $ if
$\frac{3|u|}{2}< |h|$.

If $\al(x)<1$, the situation is even simpler as we can drop $\grad
r_y^{\lam,\varepsilon}(u)\cdot h 1_{(|h|\leq 1)}$ term from
$J(u,h)$.

\ms

When $|h|>1>\frac{3|u|}{2}$, we have $|u+h|>\frac{|u|}{2}$, so
\[
r^{\lam,\varepsilon}_{y}(u+h)\leq \frac{c_{15}}{|u|^{d-\al(y)}}.
\]
Therefore
\begin{align*}
I_{32}
&\leq \frac{c_{16}}{|u|^{d-\al(y)}}\int_{|h|>1}|h|^{-d-\al(x)+\eps}dh\\
&\leq
\frac{c_{16}}{|u|^{d-\al(y)}}\int_{|h|>\frac{3|u|}{2}}|h|^{-d-\al(x)+\eps}dh.\\
&\leq \frac{c_{17}}{|u|^{d+\al(x)-\al(y)-\eps}}.
\end{align*}

Case 2: If $\frac{3|u|}{2}\ge 1$, then we have
\[
I_3 \leq
c_{18}~r^{\lam,\varepsilon}_{y}(u)\int_{|h|>\frac{3|u|}{2}}|h|^{-d-\al(x)-\eps}dh
\leq \frac{c_{19}}{|u|^{d+\al(x)-\al(y)-\eps}},
\]
since
\[
\displaystyle\sup_{\R^d\backslash
B(u,3|u|/2)}|r^{\lam,\varepsilon}_{y}(x)|\leq c~
r^{\lam,\varepsilon}_{y}(u)
\]
and $|h|^{-d-\al(x)}\leq
|h|^{-d-\al(x)-\eps}$ when $|h|>1$.

Since $|u|^{\eps}\leq |u|^{\eps/2}$ when $|u|\leq 1$, summing the
above gives
\[
|(\sL-\sM_{x})r^{\lam,\varepsilon}_{y}(u)| \leq  k_1
\frac{1}{|u|^{d+\al(x)-\al(y)-\eps/2}},\qq |u|\leq 1.
\]
We finish the proof for the first assertion of the proposition by
setting $\eta=\eps/2$. Similar arguments prove the estimate for
large $|u|$. The second assertion can be similarly proved by using
Proposition \ref{est-on-r^eps} when we estimate $|J(u,h)|$. We leave
the details to the reader. \qed

We set
\[
\sM_{z}^{y} f(x) = \int[f(x+h)-f(x)-1_{(|h|\leq 1)}\nabla f(x)\cdot
h]\frac{\xi (y)}{|h|^{d+\alpha(z)}}dh.\label{mixed_gen}
\]

\begin{proposition}\label{pertur-on-mixed}
If Assumption \ref{tempassum} holds, there exist a positive constant
$\kappa_{2} \in (0.\infty)$ such that
\begin{eqnarray*}
  |(\sM_{x}-\sM_x^y)r^{\lam,\varepsilon}_{y}(u)| &\leq & \kappa_2 \tfrac{|\xi(x)-\xi(y)|}{|u|^{d}},\qq |u|\leq 1,\\
  |(\sM_{x}-\sM_x^y)r^{\lam,\varepsilon}_{y}(u)| &\leq & \tfrac{\kappa_2}{\lam}|u|^{-d-2\underline{\al}},\qq |u|>
  1.
\end{eqnarray*}

In particular, for all $u$
\[
|(\sM_{x}-\sM_x^y)r^{\lam,\varepsilon}_{y}(u)|\leq
\tfrac{\kappa_2}{\lam}|u|^{-d-2\underline{\al}}.
\]
\end{proposition}

\proof The proof follows closely the proof of Proposition
\ref{pertur-on-lc} and the fact that $|u|^{-\al(x)}\le
|u|^{-\underline{\al}}$ when $|u|\ge 1$, where
$\underline{\al}=\displaystyle\inf_{x}\al(x)$.

\qed

Here is another estimate.

\begin{proposition}\label{pertur-on-ly}
If Assumption \ref{uniqassum} and \ref{tempassum} hold, there exist
a positive constant $\kappa_{3} \in (0.\infty)$ such that
\begin{eqnarray*}
  |(\sM_{x}^y-\sM_{y})r^{\lam,\varepsilon}_{y}(u)|\leq
\kappa_{3}\,\frac{|\al(x)-\al(y)|}{|u|^{d+|\al(x)-\al(y)|}}\big|\ln \tfrac{|u|}{2}\big|\qq |u| \leq 1 ,\\
  |(\sM_{x}^y-\sM_{y})r^{\lam,\varepsilon}_{y}(u)| \leq
\frac{\kappa_{3}}{\lam}\,\frac{1}{|u|^{d+2(\al(x)\land
\al(y))}}\big|\ln \tfrac{|u|}{2}\big| \qq |u|>1.
\end{eqnarray*}
In particular, for all $u$
\[
|(\sM_{x}^y-\sM_{y})r^{\lam,\varepsilon}_{y}(u)| \leq
\frac{\kappa_{3}}{\lam}\,\frac{1}{|u|^{d+2(\al(x)\land
\al(y))}}\big|\ln \tfrac{|u|}{2}\big|.
\]
\end{proposition}

\proof The only difference between the previous proposition and this
one is how we do the perturbation. In the previous Proposition
\ref{pertur-on-lc}, the difference of the kernels of the two
operators is
\[
\frac{n(x,h)-\xi(x)}{|h|^{d+\al(x)}}.
\]
Here the difference of kernels between two operators is
\[
\frac{\xi(y)}{|h|^{d+\al(x)}}-\frac{\xi(y)}{|h|^{d+\al(y)}}.
\]
The proof we carry out is similar to that of Proposition
\ref{pertur-on-lc}.

By Assumption \ref{tempassum}, $\xi(x)$ is bounded above.

When $|u|\leq 1$, where $J(u,h)$ is defined in (\ref{J(u,h)}), we
have
\begin{align*}
|(\sM_x^y-\sM_y)r^{\lam,\varepsilon}_{y}(u)|
&= \Big| \int J(u,h)\Big[\frac{\xi(y)}{|h|^{d+\al(x)}}-\frac{\xi(y)}{|h|^{d+\al(y)}}\Big]dh \Big|\\
& \leq c_1 \int_{|h|\leq \frac{|u|}{2}}\Big|J(u,h)[\frac{1}{|h|^{d+\al(x)}}-\frac{1}{|h|^{d+\al(y)}}]\Big|dh\\
&~~~+c_1 \int_{\frac{|u|}{2}<|h|\leq \frac{3|u|}{2}}
\Big|J(u,h)[\frac{1}{|h|^{d+\al(x)}}-\frac{1}{|h|^{d+\al(y)}}]\Big|dh\\
&~~~+c_1 \int_{|h|>\frac{3|u|}{2}}\Big|J(u,h)[\frac{1}{|h|^{d+\al(x)}}-\frac{1}{|h|^{d+\al(y)}}]\Big|dh\\
&=I_1+I_2+I_3.
\end{align*}

Without loss of generality, we may assume that $\al(x)>\al(y)$ in
the following proof.

First of all, looking at $I_1$, by Proposition \ref{est-on-r^eps} we
have
\begin{align*}
I_1 &\leq
c_{2}\sup_{B(u,|u|/2)}|\partial^{2}r^{\lam,\varepsilon}_{y}(z)|\int_{|h|\leq
\frac{|u|}{2}}\frac{1}{|h|^{d+\al(x)-2}}\Big|1-|h|^{\al(x)-\al(y)}\Big|dh\\
&\leq c_{3}\frac{|\al(x)-\al(y)|}{|u|^{d-\al(y)+2}}\int_{|r|\leq
\frac{|u|}{2}}\Big||r|^{2-\al(x)}\ln |r|\Big| dr\\
&\leq
c_{4}\frac{|\al(x)-\al(y)|}{|u|^{d+\al(x)-\al(y)}}\int_{|r|\leq
\frac{|u|}{2}}\Big|\ln |r|\Big|dr.\\
&\leq c_{5}\frac{|\al(x)-\al(y)|}{|u|^{d+\al(x)-\al(y)}},
\end{align*}
using $|r|^{2-\al(x)}\leq |u|^{2-\al(x)}$ on $r\leq \frac{|u|}{2}$
and the integrability of $\ln x$ when $x$ is small.

Next for $I_2$, by Proposition \ref{est-on-r^eps}, we have
\[
|\grad r_y^{\lam,\varepsilon}(u)\cdot h 1_{(|h|\leq 1)}|\le c_6
|u|^{-d+\al(y)},\qq \tfrac{|u|}{2}\le |h|\le \tfrac{3|u|}{2}.
\]

Therefore we have
\begin{align*}
I_2& \leq c_{7}\int_{\frac{|u|}{2}<|h|\leq
\frac{3|u|}{2}}\Big\{\frac{1}{|u+h|^{d-\al(y)}}+\frac{1}{|u|^{d-\al(y)}}\Big\}
\frac{\big|1-|h|^{\al(x)-\al(y)}\big|}{|h|^{d+\al(x)}}
dh\\
&\leq
c_{8}\frac{|\al(x)-\al(y)|}{|u|^{d+\al(x)}}\int_{\frac{|u|}{2}<|h|\leq
\frac{3|u|}{2}}\frac{1}{|u+h|^{d-\al(y)}}\big|\ln |h|\big|dh\\
&~~~~~~~+c_{9}\frac{|\al(x)-\al(y)|}{|u|^{d+\al(x)-\al(y)}}\int_{\frac{|u|}{2}<|h|\leq
\frac{3|u|}{2}}\big|\ln|h|\big|dh\\
&\leq
c_{10}\frac{|\al(x)-\al(y)|}{|u|^{d+\al(x)-\al(y)}}\big|\ln|\frac{u}{2}|\big|.
\end{align*}

Next for $I_3$, there are two cases.\\

Case 1: If $\frac{3|u|}{2}< 1$, we break up $I_3$ as follows:
\begin{align*}
I_3&=\int_{\frac{3|u|}{2}<|h|\leq
1}\Big|J(u,h)\Big[\frac{c_{1}}{|h|^{d+\al(x)}}-\frac{c_{1}}{|h|^{d+\al(y)}}\Big]\Big|dh\\
&~~~+\int_{|h|>
1}\Big|J(u,h)\Big[\frac{c_{1}}{|h|^{d+\al(x)}}-\frac{c_{1}}{|h|^{d+\al(y)}}\Big]\Big|dh\\
&=I_{31}+I_{32}.
\end{align*}
Then
\begin{align*}
I_{31}&\leq c_{11}\sup_{\R^d\backslash
B(u,\tfrac{3|u|}{2})}|\partial^{2}r^{\lam,\varepsilon}_{y}(x)|\,|\al(x)-\al(y)|\int_{\frac{3|u|}{2}<|h|\leq
  1}\Big|\frac{\ln |h|} {|h|^{d-2+\al(x)}}\Big|dh \\
  &\leq c_{12}\frac{ |\al(x)-\al(y)|}{|u|^{d+\al(x)-\al(y)}}\int_{|h|\leq
  1}\big|\ln |h|\big|dh\\
  &\leq c_{13}\frac{ |\al(x)-\al(y)|}{|u|^{d+\al(x)-\al(y)}}.
\end{align*}
When $|h|>1>\frac{3|u|}{2}$, we have $|u+h|>\frac{|u|}{2}$, so
\[
r^{\lam,\varepsilon}_{y}(u+h)\leq \frac{c_{14}}{|u|^{d-\al(y)}}.
\]
Recall that $\underline{\al}=\displaystyle\inf_{x}\al(x)$.

Then
\begin{align*}
I_{32}
&\leq c_{15}\frac{ |\al(x)-\al(y)|}{|u|^{d-\al(y)}}\int_{|h|>1}|h|^{-d-\al(x)}\big|\ln |h|\big|dh\\
&\leq
c_{16}\frac{ |\al(x)-\al(y)|}{|u|^{d-\al(y)}}\int_{|h|>\frac{3|u|}{2}}|h|^{-d-\al(x)+\underline{\al}/2}dh\\
&\leq c_{17}\frac{ |\al(x)-\al(y)|}{|u|^{d+\al(x)-\al(y)-\underline{\al}/2}}\\
&\leq c_{18}\frac{ |\al(x)-\al(y)|}{|u|^{d+\al(x)-\al(y)}},
\end{align*}
using $|u|^{\underline{\al}/2}\leq 1$ when $|u| \leq 1$.

Case 2: If $\frac{3|u|}{2}\ge 1$, then we have
\begin{align*}
I_2 &\leq
c_{19}|\al(x)-\al(y)|r^{\lam,\varepsilon}_{y}(u)\int_{|h|>\frac{3|u|}{2}}|h|^{-d-\al(x)}\big|\ln
|h|\big|dh\\
&\leq c_{20}\frac{|\al(x)-\al(y)|}{|u|^{d-\al(y)}}
\int_{|h|>\frac{3|u|}{2}}|h|^{-d-\al(x)+\underline{\al}/2}dh \\
&\leq c_{21}\frac{ |\al(x)-\al(y)|}{|u|^{d+\al(x)-\al(y)}},
\end{align*}
since $|u|^{\underline{\al}/2}\leq 1$ when $|u| \leq 1$.

Summing up the above proves the estimate for $|u|$ small. Following
similar arguments proves the estimate for large $|u|$. \qed

\begin{lemma}\label{2ndlemma}
If $r \leq 1$ and $\beta(r)$ is defined as in Assumption
\ref{uniqassum} and satisfies the condition in Assumption
\ref{uniqassum} (b), then there exists a constant $\kappa_{4}$ such
that

\begin{description}
\item(i) $r^{-1-\beta(r)+\eps} \leq \kappa_{4} r^{-1+\eps}$;
\item(ii) $\beta(r) r^{-1-\beta(r)} \leq \kappa_{4} \frac{\beta(r)}{r}$.
\end{description}
\end{lemma}
\proof By Assumption \ref{uniqassum} (c), $\beta(r)\ln(r)\rightarrow
0$ as $r\rightarrow 0$, and then $r^{\beta(r)}\rightarrow 1$ as
$r\rightarrow 0$. The lemma follows easily. \qed

\begin{proposition}\label{pertur-norm-est}
Suppose Assumption \ref{uniqassum} and \ref{tempassum} hold. Let $g
\in C^2$ with compact support. There exists a positive constant
$\widetilde{\lam}$ such that
\[
\Big|\int(\sL-\sM_y)r^{\lam,\varepsilon}_{y}(x-y)g(y)dy\Big|\leq
\frac{1}{2}\|g\|,\q \mathrm{if}\,\lam\ge\widetilde{\lam}.
\]
\end{proposition}

\proof

\begin{align*}
\Big|\int(\sL-\sM_y)r^{\lam,\varepsilon}_{y}(x-y)g(y)dy\Big|
&\le \Big|\int(\sL-\sM_x)r^{\lam,\varepsilon}_{y}(x-y)g(y)dy\Big|\\
&~~~+\Big|\int(\sM_x-\sM_x^y)r^{\lam,\varepsilon}_{y}(x-y)g(y)dy\Big|\\
&~~~+\Big|\int(\sM_x-\sM_y)r^{\lam,\varepsilon}_{y}(x-y)g(y)dy\Big|\\
&~~~=\sJ_1+\sJ_2+\sJ_3.
\end{align*}

By Proposition \ref{pertur-on-lc}, Lemma \ref{2ndlemma} and the fact
$\underline{\al}=\displaystyle\inf_x\al(x)<2$,
\begin{align*}
\sJ_1 &\le \|g\| \int_{|x-y|\lam^{\frac{1}{4}}\le 1}
\frac{\kappa_1}{|x-y|^{d+\al(x)-\al(y)-\eta}}dy\\
&~~~~~~~~~+\|g\|\int_{|x-y|\lam^{\frac{1}{4}}>
1}\frac{\kappa_1}{\lam}\frac{1}{|x-y|^{d+\al(x)+\al(y)}}dy\\
&\le
\|g\|\int_{0}^{\lam^{-\frac{1}{4}}}\frac{c_1}{r^{\beta(r)+1-\eta}}dr
+c_1\|g\| \displaystyle\int_{\lam^{-\frac{1}{4}}}^{\infty}\lam^{-1}
\frac{1}{r^{1+2\underline{\al}}}dr\\
&\le c_2 \|g\|(\lam^{-\frac{\eta}{4}}+\lam^{-1+\underline{\al}/2}).\\
\end{align*}

Taking $\lam$ large enough, say $\lam\ge \lam_1$, such that
\[
\sJ_1 \leq \frac{1}{6}\|g\|.
\]

By Proposition Assumption \ref{uniqassum} and Proposition
\ref{pertur-on-mixed},
\begin{align}
\sJ_2 &\le \|g\| \int_{|x-y|\lam^{\frac{1}{4}}\le 1}
\kappa_2\frac{\xi(x)-\xi(y)}{|x-y|^{d}}dy \nn \\
&~~~~~~~~~+\|g\|\int_{|x-y|\lam^{\frac{1}{4}}>
1}\frac{\kappa_2}{\lam}\frac{1}{|x-y|^{d+2\underline{\al}}}dy\nn\\
&\le c_3\|g\|\int_{0}^{\lam^{-\frac{1}{4}}}\frac{\psi(r)}{r}dr
+c_3\|g\| \displaystyle\int_{\lam^{-\frac{1}{4}}}^{\infty}\lam^{-1}
\frac{1}{r^{1+2\underline{\al}}}dr\nn\\
&\le c_4 \|g\|\big(\lam^{-1+\underline{\al}/2}
+\int_0^1\frac{\psi(r)}{r}1_{(0,\lam^{-1/4})}(r)dr\big)\label{J2}.
\end{align}

Letting $\lam \to \infty$, the first term of \eqref{J2} goes to $0$.
By the Dini Continuity of $\xi(x)$ and the dominated convergence
theorem, the second term of \eqref{J2} also goes to $0$.

Now take $\lam\ge\lam_2$ such that

\[
\sJ_2 \leq \frac{1}{6}\|g\|.
\]

Lastly, by Proposition \ref{pertur-on-ly} and Lemma \ref{2ndlemma},

\begin{align}
\sJ_3 &\le \|g\| \int_{|x-y|\lam^{\frac{1}{4}}\le 1}
\kappa_3\frac{\al(x)-\al(y)}{|x-y|^{d+\al(x)-\al(y)}}\big|\ln \frac{|y|}{2}\big|dy \nn \\
&~~~~~~~~~+\|g\|\int_{|x-y|\lam^{\frac{1}{4}}>
1}\frac{\kappa_3}{\lam}\frac{1}{|x-y|^{d+2\underline{\al}}}\big|\ln \frac{|y|}{2}\big|dy\nn\\
&\le
c_5\|g\|\int_{0}^{\lam^{-\frac{1}{4}}}\frac{\beta(t)}{t^{1+\beta(t)}}
|\ln t |dt +c_5\|g\|
\displaystyle\int_{\lam^{-\frac{1}{4}}}^{\infty}\lam^{-1}
\frac{1}{t^{1+2\underline{\al}}}|\ln t| dt\nn\\
&\le c_6 \big(\int_0^1
\frac{\beta(t)}{t^{1+\gamma}}1_{(0,\lam^{-1/4})}(t)dt+\int_{\lam^{-1/4}}^{\infty}\lam^{-1}
\frac{1}{t^{1+2\underline{\al}}}|\ln t|dt \big).\label{j3}
\end{align}
Since $|\ln t|\le  t^{-\gamma}$ as $t \to 0$ and $|\ln t|\le
t^{\gamma}$ as $t \to \infty$ for any $\gamma>0$.

By Assumption \ref{uniqassum} and the dominated convergence theorem,
the first term of \eqref{j3} goes to $0$ as $\lam \to \infty$.

For the second term of \eqref{j3}, we only need to take care of the
convergence of integral at $0$ and $\infty$. If we choose
$r<\min(2\underline{\al},4-2\underline{\al})$ and by the dominated
convergence theorem, the second term goes to $0$ as $\lam \to
\infty$.

Now taking $\lam \ge\lam_3$ such that
\[
\sJ_3 \leq \frac{1}{6}\|g\|.
\]

The proof is completed by taking
$\widetilde{\lam}=\max(\lam_1,\lam_2,\lam_3)$. \qed

\section{Uniqueness}

Now we are ready to show the uniqueness of the solution for the
martingale problem. Let $\P_i^x$, $i=1,2$ be two solutions to the
martingale problem starting at $x$. Let $\sR_i$ be the corresponding
resolvents.

If $f\in C^2$ with bounded first and second derivatives, by the definition of the martingale problem\\
\[
f(X_{t})-f(X_{0})-\int_{0}^{t}\sL f(X_{s})ds=\mathrm{martingale}.
\]
Taking expectations with respect to $\P_{i}^x$,\\
\begin{center}
$\E_{i}f(X_{t})-f(x)=\mathbb{E}_{i}\int_{0}^{t}\mathcal{L}f(X_{s})ds$.\\
\end{center}

Multiplying by $\lam e^{-{\lam}t}$, integrating over $t$ from
0 to $\infty$, and using Fubini's theorem gives for $i=1,2$\\
\begin{align*}
\lam\sR_{i}f-f(x)
&= \bE_{i}\int_{0}^{\infty}\lam e^{-\lam t}\int_{0}^{t}\sL f(X_{s})ds\,dt\\
&= \bE_{i}\int_{0}^{\infty}\int_{s}^{\infty}\lam e^{-\lam t}\sL f(X_{s})dt\,ds\\
&= \bE_{i}\int_{0}^{\infty}e^{-\lam s}\sL f(X_{s})ds \\
&= \sR_{i}(\sL f).
\end{align*}

Then we have

\begin{align}
\sR_{i}((\lam-\sL)f)=f(x) ~for~ i=1,2. \label{idmatrix}
\end{align}

Set $\sR_{\triangle}=\sR_{1}-\sR_{2}$. Taking the difference in
(\ref{idmatrix}) implies
\begin{align}
\sR_{\triangle}((\lam-\sL)f)=0.\label{resolventdifference}
\end{align}

Let $g$ be a $C^{2}$ function with compact support and let
\begin{align}
f_{\varepsilon}(x)=\int(r_{y}^\lambda*\varphi_{\varepsilon})(x-y)g(y)dy.\label{f-eps}
\end{align}
Note $f_{\varepsilon}(x)$ is in $C^{2}$ with bounded first and
second derivatives. Applying (\ref{resolventdifference}), we have
\begin{align}
\sR_{\triangle}((\lam-\sL)f_{\varepsilon})=0.\label{resolventdifferenceofC^2}
\end{align}

\begin{flushleft}
\emph{\textbf{Proof of Theorem \ref{main thm}.}}
\end{flushleft}
Set
\[
\th=||\sR_{\triangle}||=\display\sup_{||f||\leq
1}|\sR_{\triangle}f|.
\]
Note $|\sR_{\triangle}f|\leq (2/\lam) ||f||$, so $\th<\infty$.
\q

Let $g$ be a $C^{2}$ function with compact support and let
$f_{\varepsilon}$ be defined by (\ref{f-eps}).
From (\ref{resolventdifferenceofC^2}), we get
\begin{align*}
&~~~~|\sR_{\triangle}g|\\
&= |\sR_{\triangle}((\lam-\sL)f_{\varepsilon})-\mathcal{R}_{\triangle}g|\\
&=
|\sR_{\triangle}\int(\lam-\sM_{y}+\sM_{y}-\sL)r_{y}^{\lambda,\varepsilon}(x-y)g(y)dy-\sR_{\triangle}g|\\
&\leq  |\sR_{\triangle}(A_1-g)| +
|\sR_{\triangle}A_2|\\
&= I_1+I_2,
\end{align*}

where
\begin{eqnarray*}
  A_1 &=& \int(\lam-\sM_{y})r_{y}^{\lam,\varepsilon}(x-y)g(y)dy \\
  A_2 &=& \int(\sM_{y} -
  \sL)r_{y}^{\lam,\varepsilon}(x-y)g(y)dy.
\end{eqnarray*}

First, we look at $I_1$. Since
$(\lam-\sM_{y})r_{y}^{\lam,\varepsilon}=\varphi_{\varepsilon}$, then
\[
A_1=\int\varphi_{\varepsilon}(x-y)g(y)dy=g*\varphi_{\varepsilon}.
\]

We have
\[
\displaystyle\limsup_{\varepsilon \to
0}|\sR_{\triangle}(g*\varphi_{\varepsilon})-\sR_{\triangle}g|=0,
\]
since $g*\varphi_{\varepsilon}\to g$ uniformly.

Finally, let us look at $I_2$. By Proposition \ref{pertur-norm-est}
and take $\lam \ge \widetilde{\lam}$, we have
\[
I_2 \leq \frac{1}{2}\th \|g\|.
\]

Then
\begin{align*}
|\sR_{\triangle}g|
&\leq |\sR_{\triangle}(A_1+A_2)| \\
&\leq  \frac{1}{2}\th ||g||.
\end{align*}

Taking the sup over ${g\in C_{c}^{\infty}}$, what we get is $\th
\leq \frac{1}{2}\th $. Since $\theta<\infty$, we must have
$\theta=0$, i.e. $\sR_{1}f=\sR_{2}f$. By the uniqueness of the
Laplace transform, we have $\E_{1}f(X_t)=\E_{2}f(X_t)$ for almost
every t. Since the paths of $X_t$ are right continuous and $f$ is
continuous, then we have equality for all t. That the finite
dimensional distributions under $\P^x_1$ and $\P^x_2$ are the same
for each $x$ now follows by using the Markov property.

Lastly, we need to do a localization argument. Since $\xi(x)$ is
Dini Continuous, there must be a neighborhood of $x_0$ such that
Assumption \ref{tempassum} holds. This means that we have local
uniqueness for the martingale problem for $\sL$ started at $x_0$.
Then we follow some standard arguments; see, e.g., Chapter VI of
\cite{Ba4} to complete the proof of Theorem \ref{main thm}. \qed

\begin{remark}\label{remark:1}
{\rm We could actually compute the constant $c_{1}$ in Proposition
\ref{tranden} as a consequence of calculating $c_{d,k\al}$  by
looking at the function $f(x)=e^{-|x|^2/2}.$ Use polar coordinates
to get
\[
\omega_{d-1}(2\pi)^{d/2}\int_{0}^{\infty}r^{z+d-1}e^{-|r|^2/2}dr=c_{d,z}\omega_{d-1}\int_{0}^{\infty}r^{-z-d+d-1}e^{-|r|^2/2}dr.
\]
Do a change of variable $s=r^2/2$ and use the definition of the
gamma function to get
\[
c_{d,z}=2^{(d+z)}\pi^{d/2}\frac{\Gamma((z+d)/2)}{\Gamma(-z/2)}.
\]
Replacing $z$ by $k\al$, we have
\[
c_{d,k\al}=(2\pi)^{-d}2^{(d+k\al)}\pi^{d/2}\frac{\Gamma((d+\al)/2)}{\Gamma(-k\al/2)}=2^{k\al}\pi^{-d/2}
\frac{\Gamma((d+k\al)/2)}{\Gamma(-k\al/2)}.
\]
In particular,
\[
c_{1}=-c_{d,\al}=2^{\al}\pi^{-d/2}\frac{\Gamma((d+\al)/2)}{-\Gamma(-\al/2)}
\]
We see this makes sense since $\Gamma(-\al/2)$ is finite for
$0<\al<2$. }
\end{remark}

\noindent \textbf{Acknowledgment}. This paper is based on my Ph.D.
dissertation. I would like to express my gratitude to my adviser
Professor Richard Bass. His knowledge and guidance were invaluable
to me and greatly appreciated.

\def\cprime{$'$}

\ms

\begin{minipage}[t]{0.55\textwidth}
Department of Mathematics\\
University of Connecticut\\
Storrs, CT 06269-3009, USA\\
{\it huili@math.uconn.edu}
\end{minipage}

\end{document}